 \documentclass[11pt]{article}
 \usepackage[top=1in,bottom=1in,left=1.5in,right=1.5in]{geometry}
  \usepackage{amsmath,amssymb}
  \usepackage{latexsym}% defines $\Box$ for LaTeX2e
 \usepackage[dvips]{pstricks} % PSTricks
  \usepackage{pst-node}
  \usepackage[dvips]{graphicx}
  \newcommand{\R}{\mathbb{R}}

  \newcommand{\s}{\mathbf{s}}

  \newcommand{\U}{\mathbf{U}}
  \renewcommand{\u}{\mathbf{u}}
  
  \newcommand{\V}{\mathbf{V}}

  \newcommand{\x}{\mathbf{x}}
  
  \newcommand{\y}{\mathbf{y}}
  
  \newcommand{\z}{\mathbf{z}}
  \newcommand{\0}{\mathbf{0}}
  \newcommand{\1}{\mathbf{1}}
  
  \newcommand{\bet}{\mbox{\boldmath{$\beta$}}}
  \newcommand{\et}{\mbox{\boldmath{$\eta$}}}
  \newcommand{\bxi}{\mbox{\boldmath{$\xi$}}}

  \newcommand{\lam}{\mbox{\boldmath{$\lambda$}}}

  \newcommand{\cA}{\mathcal{A}}
  \newcommand{\cB}{\mathcal{B}}
  \newcommand{\cC}{\mathcal{C}}

  \newcommand{\lan}{\langle}
  \newcommand{\ran}{\rangle}
  \newcommand{\an}[1]{\lan#1\ran}
  
  \newcommand{\hs}{\hspace*{\parindent}}
  \newcommand{\proof}{\hs \textbf{Proof.\ }}

  \newcommand{\trans}{^\top}
  \newcommand{\qed}{\hspace*{\fill} $\Box$\\}

  \newcommand{\rH}{\mathrm{H}}

  \newtheorem{theo}{\bfseries \hs Theorem}[section]

  \newtheorem{lemma}[theo]{\bfseries \hs Lemma}
  \newtheorem{corol}[theo]{\bfseries \hs Corollary}

  \numberwithin{equation}{section} % Automatically number equations within sections

\begin{document}

 \title{Positive diagonal scaling of a nonnegative tensor \\to one with prescribed slice sums}
 \author{
 Shmuel Friedland\thanks{This research started during author's participation in AIM workshop
 ``Nonnegative Matrix Theory: Generalizations and Applications",
 December 1--5, 2008.}\\
 Department of Mathematics, Statistics and Computer Science\\
 University of Illinois at Chicago\\ Chicago, Illinois 60607-7045,
 USA\\ \texttt{E-mail: friedlan@uic.edu}
 }

 \renewcommand{\thefootnote}{\arabic{footnote}}
 \date{January 16, 2010 }
 \maketitle
 \begin{abstract}
 In this paper we give necessary and sufficient conditions on a nonnegative tensor
 to be diagonally equivalent to a tensor with prescribed slice sums.
 These conditions are variations of Bapat-Raghavan and Franklin-Lorenz conditions.
 \end{abstract}

%\mpar{check,more. I agree SF}


 \section{Introduction}\label{intro}
 For a positive integer $m$ let $\an{m}$ be the set $\{1,\ldots,m\}$.
 For positive integers $d,m_1,\ldots,m_d$
 denote by $\R^{m_1\times\ldots\times m_d}$ the linear space $d$-mode tensors
 $\cA=[a_{i_1,i_2,\ldots,i_{d}}], i_j\in \an{m_j}, j\in\an{d}$.  Note that a $1$-mode tensor is a vector, and a $2$-mode tensor is a matrix.
 Assume that $d\ge 2$.  For a fixed $i_k\in\an{m_k}$ the $(d-1)$-mode tensor $[a_{i_1,\ldots,i_d}], i_j\in\an{m_j},j\in \an{d}\backslash
 \{k\}$ is called the $(k,i_k)$ \emph{slice} of $\cA$.
 For $d=2$ the $(1,i)$ slice and the $(2,j)$ slice are the $i-th$ row and the $j-th$ column of a given matrix.
 Let
 \begin{equation}\label{kikslicesum}
 s_{k,i_k}:=\sum_{i_j\in\an{m_j},j\in \an{d}\backslash\{k\}} a_{i_1,\ldots,i_d}, \; i_k\in\an{m_k},k\in\an{d}
 \end{equation}
 be the $(k,i_k)$-slice sum.  Denote
 \begin{equation}\label{kikslsumvec}
 \s_k:=(s_{k,1}, \ldots,s_{k,m_k})\trans, \quad k\in\an{d}
 \end{equation}
 the $k$-slice vector sum.  Note that $(k,i_k)$-slice sums satisfy the compatibility conditions
 \begin{equation}\label{compslccon}
 \sum_{i_1=1}^{m_1} s_{1,i_1}=\ldots=\sum_{i_d=1}^{m_d} s_{d,i_d}.
 \end{equation}

 Two $d$-mode tensors $\cA=[a_{i_1,i_2,\ldots,i_{d}}],\cB=[b_{i_1,i_2,\ldots,i_{d}}]\in
 \R^{m_1\times\ldots\times m_d}$
 are called \emph{positive diagonally} equivalent if there exist $\x_k=(x_{k,1},\ldots,
 x_{k,m_k})\trans\in\R^{m_k}, k\in\an{d}$
 such that
 $a_{i_1,\ldots,i_d}=b_{i_1,\ldots,i_d}e^{x_{1,i_1}+\ldots+x_{d,i_d}}$ for all $i_j\in\an{m_j}$
 and $j\in\an{d}$.
 Denote by  $\R_+^{m_1\times\ldots\times m_d}$ the cone of nonnegative,(entrywise),  $d$-mode tensors.

 In this paper we assume that $\cB=[b_{i_1,i_2,\ldots,i_{d}}]\in
 \R_+^{m_1\times \ldots\times m_d}$ is a given nonnegative tensor with no zero slice $(k,i_k)$.
 Let $\s_k\in \R_+^{m_k}, k\in\an{d}$ are given $k$ positive vectors satisfying the conditions
 (\ref{compslccon}).
 Denote by $\R_+^{m_1\times \ldots\times m_d}(\cB, \s_1,\ldots,\s_d)$ the set of all nonnegative
 $\cA=[a_{i_1,i_2,\ldots,i_{d}}]\in
 \R_+^{m_1\times \ldots m_d}$ having the same zero pattern as $\cB$, i.e. $a_{i_1,\ldots,i_d}=0
 \iff b_{i_1,\ldots,i_d}=0$
 for all indices $i_1,\ldots,i_d$, and satisfying the condition (\ref{kikslicesum}).
 The aim of this paper is to give new necessary and sufficient conditions on $\cB$ so that
 $\R_+^{m_1\times \ldots m_d}(\cB, \s_1,\ldots,\s_d)$
 contains a tensor $\cA$, which is positively diagonally equivalent to $\cB$.  For matrices, i.e. $d=2$,
 this problem was solved by
 Menon \cite{Men68} and Brualdi \cite{Bru68}.  See also \cite{MS69}.  For the special case of positive
 diagonal equivalence
 to doubly stochastic matrices see \cite{BPS66} and \cite{SK}.
 The result of Menon was extended for tensors independently by Bapat-Raghavan \cite{BR89} and
 Franklin-Lorenz \cite{FL89}.  (See \cite{Bap82} and \cite{Rag84} for the special case where all the
 entries of $\cB$ are positive.)\footnote{I thank Ravi Bapat for pointing out these references in
 December 2009.}  In this paper we give a different necessary and sufficient
 conditions for the solution of this problem.
 \begin{theo}\label{maintheo}  Let $\cB=[b_{i_1,i_2,\ldots,i_{d}}]\in\R_+^{m_1\times\ldots\times m_d}$,
 $(d\ge 2)$, be a
 given nonnegative tensor with no $(k,i_k)$-zero slice.  Let $\s_k\in\R_+^{m_k}, k=1,\ldots,d$ be given
 positive vectors satisfying (\ref{compslccon}).
 Then there exists a nonnegative tensor $\cA\in \R_+^{m_1\times\ldots\times m_d}$, which is positive
 diagonally equivalent
 to $\cB$ and having each $(k,i_k)$-slice sum equal to $s_{k,i_k}$, if and only the following conditions.
 The system of the inequalities and equalities for $\x_k=(x_{k,1},\ldots,\x_{k,m_k})
 \trans\in\R^{m_k}, k=1,\ldots,d$,
 \begin{eqnarray}\label{ineqcond}
 x_{1,i_1}+x_{2,i_2}+\ldots+x_{d,i_d}\le 0 \textrm{ if } b_{i_1,i_2,\ldots,i_d}>0,\\
 \label{eqcond}
 \s_k\trans \x_k=0 \textrm{ for } k=1,\ldots,d,
 \end{eqnarray}
 imply one of the following equivalent conditions
 \begin{enumerate}
 \item\label{eqcond1}  $x_{1,i_1}+x_{2,i_2}+\ldots+x_{d,i_d}= 0$ if $b_{i_1,i_2,\ldots,i_d}>0$.
 \item\label{eqcond2}
 $\sum_{b_{i_1,i_2,\ldots,i_d}>0} x_{1,i_1}+x_{2,i_2}+\ldots+x_{d,i_d}=0$.

 \end{enumerate}
 In particular, there exists at most one tensor $\cA\in \R_+^{m_1\times\ldots\times m_d}$ with $(k,i_k)$-slice
 sum $s_{k,i_k}$ for all $k,i_k$,  which is positive diagonally equivalent to $\cB$.
 \end{theo}

 The above yields the following corollary.
 \begin{corol}\label{mentens}  Let $\cB\in\R_+^{m_1\times\ldots\times m_d}$, $(d\ge 2)$, be a
 given nonnegative tensor with no $(k,i_k)$-zero slice.  Let $\s_k\in\R_+^{m_k}, k=1,\ldots,d$ be given
 positive vectors.  Then there exists a nonnegative tensor $\cC\in \R_+^{m_1\times\ldots\times m_d}$,
 which is positive diagonally equivalent
 to $\cB$ and each $(k,i_k)$-sum slice equal to $s_{k,i_k}$, if and only if there exists a nonnegative tensor
 $\cA=[a_{i_1,i_2,\ldots,i_{d}}]\in\R_+^{m_1\times\ldots\times m_d}$,
 having the same zero pattern as $\cB$, which satisfies (\ref{kikslicesum}).
 \end{corol}
 For matrices, i.e. $d=2$, the above corollary is due Menon \cite{Men68}.
 For $d= 3$ this result is due to \cite[Thm 3]{BR89} and for $d\ge 3$ \cite{FL89}.
 Brualdi in \cite{Bru68} gave a nice and simple characterization for the set of nonnegative matrices, with
 prescribed zero pattern and with given positive row and column sums, to be not empty.
 It is an open problem to find an analog of Brualdi's results for $d$-mode tensors, where $d\ge 3$.

 Note that the conditions of Theorem \ref{maintheo}  are stated as a
 linear programming problem.  Hence the existence of a positive diagonally equivalent tensor $\cA$ can
 be determined in polynomial time.  If such $\cA$ exists,  we show that it can be found
 by computing the unique minimal point of certain strictly convex functions $f$.
 Hence, Newton method can be applied to find the unique minimal point of $f$ and its value very fast.
 (See \S3.)

 \section{Proof of the main theorem}

 Identify  $\R^{m_1}\times\R^{m_2}\times\ldots\times\R^{m_d}$ with $\R^n$, where $n=\sum_{k=1}^d m_k$.
 We view $\y\in\R^n$ as a vector $(\x_1\trans,\ldots,\x_d\trans)\trans$, where $\x_k\in\R^{m_k}, k\in\an{d}$.
 Let $\|\y\|:=\sqrt{\y\trans\y}$.  Define
 \begin{equation}\label{deffunct}
 f(\y)=f((\x_1\trans,\ldots,\x_d\trans)\trans):=\sum_{i_j\in\an{m_j},j\in\an{d}}
 b_{i_1,\ldots,i_d}e^{x_{1,i_1}+\ldots+x_{d,i_d}}.
 \end{equation}
 Clearly, $f$ is a convex function on $\R^n$.
 Denote by $\U(\s_1,\ldots,\s_d)\subset \R^n$ the subspace of vectors  $(\x_1\trans,\ldots,\x_d\trans)\trans$
 satisfying the equalities (\ref{eqcond}).
 \begin{lemma}\label{critptf}  Let $\cB=[b_{i_1,i_2,\ldots,i_{d}}]\in\R_+^{m_1\times\ldots\times m_d}$,
 $(d\ge 2)$, be a given nonnegative tensor with no $(k,i_k)$-zero slice.  Let $\s_k\in\R_+^{m_k}, k=1,\ldots,d$
 be given positive vectors satisfying (\ref{compslccon}).  Then there exists a nonnegative tensor
 $\cA\in \R_+^{m_1\times\ldots\times m_d}$,
 which is positive diagonally equivalent to $\cB$ and having each $(k,i_k)$-slice sum equal to $s_{k,i_k}$, if
 and only the restriction of $f$ to the subspace $\U(\s_1,\ldots,\s_d)$, ($f|\U(\s_1,\ldots,\s_d)$, has a
 critical point.

 \end{lemma}
 \proof  Assume first that $f|\U(\s_1,\ldots,\s_d)$ has a critical point.  Use Lagrange multipliers, i.e.
 consider the function $f-\sum_{k=1}^d \s_k\trans\x_k$, to deduce
 the existence of $\lam=(\lambda_1,\ldots,\lambda_d)\trans$ and $(\bxi_1\trans,\ldots,\bxi_d\trans)
 \trans\in\U(\s_1,\ldots,\s_d)$,
 where $\bxi_k=(\xi_{k,1},\ldots,\xi_{k,i_k})\trans, k\in\an{d}$, satisfying the following conditions
 \begin{equation}\label{critcon}
 \sum_{i_j\in\an{m_j},j\in \an{d}\backslash\{k\}} b_{i_1,\ldots,i_d} e^{\xi_{1,i_1}+\ldots+\xi_{d,i_d}}=
 \lambda_k s_{k,i_k},\;, \; i_k\in\an{m_k},k\in\an{d}.
 \end{equation}
 Since $\s_k>$ is a positive vector and the $(k,i_k)$-slice of $\cB$ is not a zero slice we deduce that
 $\lambda_k>0$.
 Summing up the above equation on $i_k=1,\ldots,m_k$, and using the equalities (\ref{compslccon})
 we deduce that $\lambda_1=\ldots=\lambda_d>0$.  Then $\cA=[b_{i_1,\ldots,i_d}e^{(\xi_{1,i_1}-\log\lambda_1)+
 \xi_{2,i_2}+\ldots+\xi_{d,i_d}}]$.

 Vice versa suppose $\cA=[b_{i_1,\ldots,i_d}e^{x_{1,i_1}+\ldots+x_{d,i_d}}]$ has $(k,i_k)$-slice sum equal to
 $s_{k,i_k}$ for all $(k,i_k)$.  Let $\1_m=(1,\ldots,1)\trans\in\R^m$.  Then there exists a unique
 $t_i\in\R$ such that $\s_k\trans(\x_k-t_k\1_{m_k})=0$
 for $k\in\an{d}$.  Let $\bxi_k:=\x_k-t_k\s_k, k\in\an{d}$.  Then (\ref{critcon}) holds.  \qed

 Denote by $\V(\s_1,\ldots,\s_d)\subset \U(\s_1,\ldots,\s_d)$ the subspace of all vectors
 $(\x_1\trans,\ldots,\x_d\trans)\trans$
 satisfying the condition \emph{\ref{eqcond1}} of Theorem \ref{maintheo}.
 Clearly, for each $\y\in\R^n$ the function $f$ has a constant value $f(\y)$ on the affine set $\y+
 \V(\s_1,\ldots,\s_d)$.  Hence, if $\et\in\U(\s_1,\ldots,\s_d)$ is a critical point of $f|\U(\s_1,\ldots,\s_d)$
 then any point in $\et+\V(\s_1,\ldots,\s_d)$ is also a critical of  $f|\U(\s_1,\ldots,\s_d)$.
 Denote by $\V(\s_1,\ldots,\s_d)^{\perp}\subset \U(\s_1,\ldots,\s_d)$, the orthogonal complement of
 $\V(\s_1,\ldots,\s_d)$ in $\U(\s_1,\ldots,\s_d)$.  Thus, to study the existence of the critical points of
 $f|\U(\s_1,\ldots,\s_d)$, it is enough
 to study the existence of the critical points of $f|\V(\s_1,\ldots,\s_d)^{\perp}$.
 Since the function $e^{at}$ is strictly
 convex for $t\in\R$ for any $a\ne 0$, more precisely $(e^{at})''=a^2e^{at}>0$  we deduce the following.
 \begin{lemma}\label{strictconvf}  Let $\cB=[b_{i_1,i_2,\ldots,i_{d}}]\in\R_+^{m_1\times\ldots\times m_d}$,
 $(d\ge 2)$, be a given nonnegative tensor with no $(k,i_k)$-zero slice.  Let $\s_k\in\R_+^{m_k}, k=1,\ldots,d$
 be given positive vectors satisfying (\ref{compslccon}).  Let $\U(\s_1,\ldots,\s_d),\V(\s_1,\ldots,\s_d),
 \V(\s_1,\ldots,\s_d)^{\perp}$ be defined as above.
 Then $f|\V(\s_1,\ldots,\s_d)^{\perp}$ is strictly convex.  More precisely, the Hessian matrix of
 $f|\V(\s_1,\ldots,\s_d)^{\perp}$
 has positive eigenvalues at each point of  $\V(\s_1,\ldots,\s_d)^{\perp}$.
 \end{lemma}
 \begin{theo}\label{eqivcondscal}  Let $\cB=[b_{i_1,i_2,\ldots,i_{d}}]\in\R_+^{m_1\times\ldots\times m_d}$,
 $(d\ge 2)$, be a given nonnegative tensor with no $(k,i_k)$-zero slice.  Let $\s_k\in\R_+^{m_k}, k=1,\ldots,d$
 be given positive vectors satisfying (\ref{compslccon}).  Then the following conditions are equivalent.
 \begin{enumerate}
 \item\label{eqivcondscal1}  $f|\U(\s_1,\ldots,\s_d)$ has a global minimum.
 \item\label{eqivcondscal2}  $f|\U(\s_1,\ldots,\s_d)$ has a critical point.
 \item\label{eqivcondscal3}  $\lim f(\y_l)=\infty$ for any sequence $\y_l\in\V(\s_1,\ldots,\s_k)^{\perp}$
 such that $\lim\|\y_l\|= \infty$.
 \item\label{eqivcondscal4}  The only $\y=(\x_1\trans,\ldots,\x_d\trans)\trans\in\V(\s_1,\ldots,\s_k)^{\perp}$
 that satisfies (\ref{ineqcond}) is $\y=\0_n$.

 \end{enumerate}
 \end{theo}
 \proof
 \ref{eqivcondscal1} $\Rightarrow$ \ref{eqivcondscal2}.  Trivial.

 \noindent
 \ref{eqivcondscal2} $\Rightarrow$ \ref{eqivcondscal3}.  Let $\bet\in \U(\s_1,\ldots,\s_d)$ be a critical point
 of $f|\U(\s_1,\ldots,\s_d)$.  Hence any point in $\bet+\V(\s_1,\ldots,\s_d)$ is a critical point of
 $f|\U(\s_1,\ldots,\s_d)$.  Hence $f|\U(\s_1,\ldots,\s_d)$ has a critical point
 $\bxi\in\V(\s_1,\ldots,\s_d)^{\perp}$.
 In particular, $\bxi$ is a critical point of $f|\V(\s_1,\ldots,\s_d)^{\perp}$.
 Let $\z\in \V(\s_1,\ldots,\s_d)^{\perp},
 \|\z\|=1$.   For $t\in\R$ define $g_{\z}(t):= f(\bxi+t\z)$.  So $g_{\z}$ is strictly convex on $\R$ and
 $g_{\z}'(0)=0$.  Let  $\rH(f)(\y)$ be the Hessian matrix of $f|\V(\s_1,\ldots,\s_d)^{\perp}$ at
 $\y\in \V(\s_1,\ldots,\s_d)^{\perp}$, i.e. the symmetric matrix of the second derivatives of
 $f|\V(\s_1,\ldots,\s_d)^{\perp}$ at $\y\in \V(\s_1,\ldots,\s_d)^{\perp}$.
 Lemma \ref{strictconvf} implies that the smallest eigenvalue $\alpha(\y)$ of $\rH(f)(\y)$ is positive.
 Clearly, $\rH(\y)$ and hence $\alpha(\y)$ are continuous on $\V(\s_1,\ldots,\s_d)^{\perp}$.  Hence
 $\min_{\|\y-\bxi\|\le 1} \alpha(\y)=2a>0$.  Therefore, $g_{\z}''(t)\ge a$ for $t\in [-1,1]$.  In particular
 $g_{\z}'(t)\ge 2at$ and $g_{\z}(t)\ge f(\bxi)+at^2$
 for any $t\in [0,1]$.  So $g_{\z}(1)\ge f(\bxi)+a$.  Since $g_{\z}'(t)$ increases on $\R$ it follows that
 $g_{\z}'(t)\ge 2a$ for $t\ge 1$.  Hence $g_{\z}(t)\ge f(\bxi)+ a+ 2a(t-1)=f(\bxi)+a(2t-1)$ for $t\ge 1$.  Thus
 $f(\bxi+\u)\ge a(2\|u\|-1)$ for any $\u\in\V(\s_1,\ldots,
 \s_d), \|\u\|\ge 1$.  Hence \ref{eqivcondscal3} holds.

 \noindent
 \ref{eqivcondscal3} $\Rightarrow$ \ref{eqivcondscal1}.  Since $f=\infty$ on
 $\partial \V(\s_1,\ldots,\s_d)^{\perp}$
 it follows that $f|\V(\s_1,\ldots,\s_d)^{\perp}$ achieves its minimum at
 $\bxi\in  \V(\s_1,\ldots,\s_d)^{\perp}$.  Clearly, for any point $\y\in \U(\s_1,\ldots,\s_d)$ there exits
 $\z\in \V(\s_1,\ldots,\s_d)^{\perp}$ such that
 $\y\in \z+\V(\s_1,\ldots,\s_d)$.  Recall that $f(\y)=f(\z)\ge f(\bxi)$.  Hence  $f(\bxi)$ is the minimum of
 $f|\U(\s_1,\ldots,\s_d)$.

 \noindent
 \ref{eqivcondscal3} $\Rightarrow$ \ref{eqivcondscal4}.  Assume to the contrary that there exists
 $\0\ne\y=(\x_1\trans,\ldots,\x_d\trans)\trans\in \V(\s_1,\ldots,\s_d)^{\perp}$
 which satisfies (\ref{ineqcond}). Hence, there exists $i_j\in \an{m_j}, j\in\an{d}$ such that
 $b_{i_1,\ldots,i_d}>0$ and $x_{1,i_1}+\ldots+x_{d,i_d}<0$.  Thus, there exist $\alpha_1,\ldots,\alpha_p<0$
 and $\beta_1,\ldots,\beta_l>0$
 such that  $f(t\y)=\gamma+\sum_{l=1}^p \beta_le^{t\alpha_l}$.
 (Each $\alpha_q$ is equal to some $x_{1,i_1}+\ldots+x_{d,i_d}<0$, where $b_{i_1,\ldots,i_d}>0$,
 and each $\beta_q$ is a sum of corresponding $b_{i_1,\ldots,i_d}>0$.)
 Hence, $\lim_{t\to\infty} f(t\y)=\gamma$, which contradicts
 \ref{eqivcondscal3}.

 \noindent
 \ref{eqivcondscal4} $\Rightarrow$ \ref{eqivcondscal3}.  Let $\y\in(\x_1\trans,\ldots,\x_d\trans)\trans\in
 \V(\s_1,\ldots,\s_d)^{\perp}, \|\y\|=1$.  Then
 $$h(\y):=\max_{b_{i_1,\ldots,i_d}>0} x_{1,i_1}+\ldots+x_{d,i_d} >0.$$  The continuity of $h(\y)$
 on the unit sphere in $\V(\s_1,\ldots,\s_d)^{\perp}$ implies that
 $$\min_{\y\in\V(\s_1,\ldots,\s_d)^{\perp}, \|\y\|=1} h(\y)=\alpha>0.$$
 Let $\beta=\min_{b_{i_1,\ldots,i_d}>0} b_{i_1,\ldots,i_d}>0$.  Hence, for any
 $\y\in(\x_1\trans,\ldots,\x_d\trans)\trans\in \V(\s_1,\ldots,\s_d)^{\perp}$, $\|\y\|=1$
 and $t>0$ we have that $f(t\y)\ge\beta e^{\alpha t}$.
 This inequality yields \ref{eqivcondscal3}.  \qed

 \textbf{Proof of Theorem \ref{maintheo}.}  Assume first that there exists a nonnegative tensor
 $\cA\in \R_+^{m_1\times\ldots\times m_d}$, which is positive diagonally equivalent
 to $\cB$ and having each $(k,i_k)$-slice sum equal to $s_{k,i_k}$.  Lemma \ref{critptf} yields that
 $f|\U(\s_1,\ldots,\s_d)$ has a critical point, i.e. the condition \emph{\ref{eqivcondscal2}}  of Theorem
 \ref{eqivcondscal} holds.
 Since $f|\V(\s_1,\ldots,\s_d)^{\perp}$ is strictly convex, it has a unique critical point
 $\bxi\in\V(\s_1,\ldots,\s_d)^{\perp}$.
 Hence all critical points of a convex $f|\U(\s_1,\ldots,\s_k)$ must be of the form $\bxi+\V(\s_1,\ldots,\s_d)$.
 The proof of Lemma \ref{critptf} yields that $\cA$ is unique.

 Theorem \ref{eqivcondscal} implies the condition \emph{\ref{eqivcondscal4}}.  Hence the conditions
 (\ref{ineqcond}) and (\ref{eqcond}) yield the conditions \emph{\ref{eqcond1}} and \emph{\ref{eqcond2}} of
 Theorem  \ref{maintheo}.

 Assume that the conditions (\ref{ineqcond}) and (\ref{eqcond}) hold.  Clearly the conditions
 \emph{\ref{eqcond1}} and \emph{\ref{eqcond2}} of Theorem  \ref{maintheo} are equivalent.  Suppose now that
 the conditions (\ref{ineqcond}) and (\ref{eqcond})
 imply the condition \emph{\ref{eqcond1}} of Theorem  \ref{maintheo}.  Hence the condition
 \emph{\ref{eqivcondscal4}} of Theorem \ref{eqivcondscal} holds.  Use the  the condition
 \emph{\ref{eqivcondscal2}} of Theorem \ref{eqivcondscal} and Lemma
 \ref{critptf} to deduce the existence of a nonnegative tensor
 $\cA\in \R_+^{m_1\times\ldots\times m_d}$, which is positive diagonally equivalent
 to $\cB$ and having each $(k,i_k)$-slice sum equal to $s_{k,i_k}$.

 \qed

  \textbf{Proof of Corollary \ref{mentens}.}  We prove the nontrivial part of the corollary.
  Suppose that there exists a nonnegative tensor
 $\cA=[a_{i_1,i_2,\ldots,i_{d}}]\in\R_+^{m_1\times\ldots\times m_d}$,
 having the same zero pattern as $\cB$, which satisfies (\ref{kikslicesum}).  Clearly, $\cA$ is positively
 diagonally equivalent to $\cA$ and has each $(k,i_k)$-sum slice equal to $s_{k,i_k}$.  Apply Theorem
 \ref{maintheo} to $\cA$ to deduce that
 the set of inequalities  $x_{1,i_1}+x_{2,i_2}+\ldots+x_{d,i_d}\le 0$ if  $a_{i_1,i_2,\ldots,i_d}>0$, together
 with the equalities (\ref{eqcond}) yields the condition
 $\sum_{a_{i_1,i_2,\ldots,i_d}>0} x_{1,i_1}+x_{2,i_2}+\ldots+x_{d,i_d}=0$.
 Since $a_{i_1,\ldots,i_d}>0 \iff b_{i_1,\ldots,i_d}>0$ we deduce that the conditions (\ref{ineqcond}) and
 (\ref{eqcond})
 of Theorem \ref{maintheo} yield the condition  \emph{\ref{eqcond2}} of Theorem  \ref{maintheo}.
 Hence there exists a nonnegative tensor $\cC\in \R_+^{m_1\times\ldots\times m_d}$, which is positive diagonally equivalent
 to $\cB$ and has $(k,i_k)$-sum slices equal to $s_{k,i_k}$.  \qed

 \section{Remarks}
 Theorem \ref{maintheo}, the main result of this paper, is stated stated in terms of linear programming.
 Hence by the results of \cite{Kha,Kar}
 one can verify these conditions in polynomial time.  The proof of Theorem \ref{eqivcondscal},
 combined Lemma \ref{critptf}, shows that to find $\cA$, which is diagonally equivalent to $\cB$, we
 need to find the minimum of the strict convex function $f|\V(\s_1,\ldots,\s_d)^{\perp}$.
 There are many numerical methods to to find the unique minimum, e.g. \cite{BV}.
 Since the Hessian at the critical point of our strict convex function has positive eigenvalues, one should
 use the Newton algorithm, or its variant as Armijo rule \cite{NW99}, to obtain the quadratic
 convergence.

 In the special case of diagonal equivalence to doubly stochastic matrices, one can performs the
 Sinkhorn scaling algorithm \cite{SK,FL89}, which converges linearly.
 It seems to the author, that even in the case of matrices, a variant of the Newton algorithm should outperform
 the Sinkhorn scaling algorithm.  The numerical aspects of the comparison between the two algorithms
 will be done somewhere else.

 \bibliographystyle{plain}

 \emph{Acknowledgement}:  I thank  St\'ephane Gaubert and Leonid Gurvits
 for useful remarks.

\end{document}